\documentstyle[psfig]{article}

\def\0{{\bf 0}}

\def\bt{\begin{theorem}}
\def\et{\end{theorem}}
\def\bp{\begin{proposition}}
\def\ep{\end{proposition}}
\def\bl{\begin{lemma}}
\def\el{\end{lemma}}
\def\bi{\begin{itemize}}
\def\ei{\end{itemize}}
\def\bd{\begin{description}}
\def\ed{\end{description}}
\def\br{\begin{remark}}
\def\er{\end{remark}}
\def\be{\begin{equation}}
\def\ee{\end{equation}}
\def\bc{\begin{corollary}}
\def\ec{\end{corollary}}
\def\bex{\begin{example}}
\def\eex{\end{example}}

\input mssymb12.tex

\newcommand{\R}{{\Bbb R}}

\def\NABLA#1{{\mathop{\nabla\kern-.5ex\lower1ex\hbox{$#1$}}}}
\def\Nabla#1{\nabla\kern-.5ex{}_#1}

\newcommand{\p}{{\partial}}

\newcommand{\qed}{\hfill$\Box$\medskip}
\newtheorem{theorem}{Theorem} 
\newtheorem{corollary}[theorem]{Corollary}
\newtheorem{lemma}[theorem]{Lemma}
\newtheorem{proposition}[theorem]{Proposition}

\newtheorem{remark}[theorem]{Remark}
\newtheorem{example}[theorem]{Example}

\title{Symplectic Rigidity for Anosov Hypersurfaces}

\author{D. Burns \and R. Hind}

\begin{document}

\maketitle

\section{Introduction}

There is a canonical exact symplectic structure on the unit tangent bundle
of a Riemannian manifold given by pulling-back the symplectic
structure on the cotangent bundle using the Riemannian metric.

The following result is a straightforward application of the
symplectic homology theory, see \cite{cfh}, and a theorem of
J. Otal \cite{otal} and C. Croke
\cite{croke}.

\begin{theorem}
If the interiors of the unit tangent bundles of two compact Riemann surfaces
of strictly negative curvature are exact symplectomorphic then the underlying
Riemann surfaces are isometric.
\end{theorem}

An exact symplectomorphism is one for which the pull-back of the canonical
Liouville form differs from itself by an exact $1$-form.

Although in this paper we will mainly be considering tangent bundles
of surfaces, the above theorem does have some generalizations to
higher dimensions. For example, using a result of U. Hamenst\"{a}dt
in \cite{ham} we get

{\bf Theorem 1'} {\it Let $M$ and $N$ be closed, strictly negatively curved
manifolds and suppose that $N$ is a locally symmetric space. If the
interiors of the unit tangent bundles of $M$ and $N$ are exact symplectomorphic
then $M$ and $N$ are isometric.}

For completeness we explain how to obtain theorems $1$ and $1'$ from
the symplectic homology theory together with the results on length
spectrum rigidity at the end of the introduction. In the meantime
we will concentrate on surfaces.

\begin{corollary}
If the interiors of the unit tangent bundles of two Riemann surfaces
of strictly negative curvature are exact symplectomorphic then the closed
symplectic manifolds are also symplectomorphic.
\end{corollary}

Although it is a weaker result than the theorem, this corollary is
still very interesting, especially in the light of work of
Y. Eliashberg and H. Hofer, see \cite{eh}, showing that there exist
$C^\infty$-small perturbations of the standard unit ball in
$\R^{2n}$ whose interiors are symplectomorphic but whose boundaries
are not.

The main purpose of this paper is to extend the above symplectic
rigidity result to a larger class of symplectic manifolds.

Fixing our underlying smooth surface $M$, the unit tangent bundles
corresponding to different metrics can symplectically be thought
of as domains in $T^* M$ with the restricted canonical symplectic
form by applying the Legendre transform. We will look at the class
of symplectic manifolds obtained by deforming in
$T^* M$ the domains corresponding to negatively curved metrics
and again restricting the canonical symplectic form.
This class of symplectic manifolds is the same as the one obtained
by fixing a unit tangent bundle and deforming the primitive of
the symplectic form in the same cohomology class.

{\bf Remark} This class of symplectic manifolds is in a sense open and
certainly contains domains in $T^* M$ which are not the Legendre
transform of Riemannian unit tangent bundles. It also contains domains
which are not symplectomorphic to Riemannian unit tangent bundles. For
example, a negatively curved Riemannian metric can be deformed to give
a non-symmetric Finsler metric with the property that a closed
geodesic in a certain free homotopy class $\gamma$ is of different
length to the unique closed geodesic in the class $-\gamma$. Now, when
one applies the Legendre transform to a unit (Riemannian or Finsler)
tangent bundle, the resulting domain in $T^* M$ has a canonical Reeb
flow on its boundary given by restricting the Liouville form to a
contact form on the boundary. The geodesic flow on the tangent bundle
(restricted to a fixed energy level) corresponds under the Legendre
transform to this Reeb flow and in particular the lengths of closed
geodesics correspond to the periods of closed orbits of the Reeb
flow. By the symplectic homology theory, these periods are invariants
of the exact symplectomorphism type of the domain.  But for a
Riemannian metric the unique closed geodesic in the opposite class to
a given closed geodesic is just the original closed geodesic traversed
in the opposite direction and has the same length. Thus the domain
obtained from our Finsler metric cannot be 
exact symplectomorphic to a Riemannian
domain. We do not know, on the other hand, whether there exist
symmetric domains in $T^* M$ which are not symplectomorphic to Riemannian
domains.

\vskip 5pt

Now, provided that a deformed domain still intersects each fiber in a
convex set one can apply the inverse Legendre transform and associate
to the domain a Finsler metric. One might hope to generalise Theorem
$1$ and say that if two domains are symplectic then their associated
Finsler metrics must be isometric.  Unfortunately though, it is easy
to construct examples showing this to be false. In particular, it is
possible for a Finsler metric to have the same length spectrum as a
Riemannian metric while still not being Riemannian itself. To see
this, we start with a Riemannian domain $W$ in $T^* M$, say
corresponding to a metric $g$ on $M$. Let $H$ be a (Hamiltonian)
function on $T^* M$ supported in a sufficiently small neighbourhood of some point $x\in \p W$.  We assume that
the induced Hamiltonian diffeomorphism $\phi$, the time-$1$ flow of
the Hamiltonian vector field, does not preserve $\p W$, and then study
the domain $W'=\phi (W)$.  Applying the inverse Legendre transform to
$W'$, provided that $H$ was sufficiently small, we will get a domain
in $TM$ which is the unit tangent bundle of a certain Finsler
metric. We want to observe that this Finsler metric is not Riemannian.
But for a Riemannian metric, the unit circle in each tangent space
$T_p M$ is an ellipse. In this case though, there are some tangent
spaces where the unit circle of the Finsler metric coincides with that
for the Riemannian metric $g$ except in a neighbourhood of some
point. Since ellipses which coincide on open sets are actually equal,
the Finsler unit circle cannot be an ellipse and hence the metric is
not Riemannian. Choosing the Hamiltonian $H$ to be supported in a
neighbourhood of two points, the above construction can be carried out
symmetrically and gives examples of symmetric Finsler domains which
are symplectomorphic to Riemannian domains but are not Riemannian.

Given the above, we will therefore seek to generalise Corollary $2$.
First we will be more specific about exactly which class of
symplectic manifolds will be considered.

Let $\mathcal{M}$ be the class of domains $W$ in $T^* M$ with
smooth boundary such that the canonical Liouville form $\lambda$
on $T^* M$ restricts to a contact form on $\p W$ whose Reeb
vector field $X$, uniquely defined by $X\rfloor d\lambda =0$ and
$\lambda (X)=1$, generates an Anosov flow on $\p W$.
We will also assume that the zero-section in $T^* M$ lies inside
$W$, and that the fibers in $T^* M$ are star-shaped.

It follows from Anosov's structural stability theorem, see \cite{anosov},
that $\mathcal{M}$ is an open set (with a topology
of smooth convergence). Furthermore, the geodesic flow of a negatively
curved metric restricts to an Anosov flow on constant energy surfaces
and so all the deformed domains described above lie in $\mathcal{M}$,
and in fact they lie in the
same connected component $\mathcal{M}^{\circ}$. We can now
generalise the above corollary as follows.

\begin{theorem}
Suppose that the interiors of two domains $W_1$ and $W_2$ in
$\mathcal{M}^{\circ}$ are exact symplectomorphic. Then the closed
symplectic manifolds are symplectomorphic and in fact the
symplectomorphism can be taken to be the restriction of a smooth
Hamiltonian diffeomorphism on $T^* M$, perhaps composed with
the differential of a diffeomorphism of $M$.
\end{theorem}

We emphasize the smoothness here as this relies on a result
in dynamical systems due to R. de la Llave and R. Moriyon
in \cite{mor2}, see also \cite{fo}.
 
Actually, the only infomation needed about the interior
symplectomorphism above is that it preserves the length
spectrum of the Reeb flow on the boundary, that is, 
corresponding closed orbits have the same length.
Therefore we also have the following result.
Let $\p \mathcal{M} =\{\p W|W\in\mathcal{M}^{\circ}\}$.

\begin{theorem}
Suppose that two hypersurfaces $\Sigma_0$ and $\Sigma_1$ in
$\p \mathcal{M}$ have the same marked length spectrum. 
Then they are connected by a smooth $1$-parameter family
of hypersurfaces in $\p \mathcal{M}$ whose Reeb flows are
all smoothly time-preserving conjugate to the flow on
$\Sigma_0$ by a smooth family of Hamiltonian diffeomorphisms.
\end{theorem}

As both hypersurfaces lie in $\p \mathcal{M}$, each free homotopy
class in $\pi_1(M)$ corresponds to the projection of a unique
closed orbit on $\Sigma_i$. The marked length spectrum is the
map associating the length of the closed orbit to the homotopy
class.

{\bf Proof of Theorem $1$}

By \cite{croke} or \cite{otal}, it suffices to show that the surfaces
have the same marked length spectrum (up to diffeomorphism).

But suppose that the interiors are symplectomorphic via an exact
symplectomorphism $\phi$. Then $\phi$ induces an isomorphism of
the first homotopy groups of the underlying surface. By a theorem
of Nielsen, see \cite{nielsen}, any such isomorphism can be
induced by a diffeomorphism of the surface. Hence we may assume
that the map $\phi _*$ on $\pi _1(M)$ is the identity.

Now we apply the symplectic homology theory, see for example
\cite{cfh}, to deduce that each (oriented) closed geodesic in the
first surface, which gives a closed orbit of the characteristic flow
on the boundary of the symplectic manifold, must correspond to a
closed orbit of the characteristic flow and hence a closed geodesic in
the second surface of the same length. In fact we are using a version
of symplectic homology which only considers closed orbits in a fixed
homology class. With this restriction, exactly the same construction
applies.  As $\phi_* =\imath$, the second geodesic is in the same free
homotopy class and so the surfaces have the same marked length
spectrums as required.

{\bf Proof of Theorem $1'$}

Again, now according to Hemenst\"{a}dt, see \cite{ham}, it suffices
to show that the manifolds $M$ and $N$ have the same marked length spectrum.
As $M$ and $N$ are strictly negatively curved, each free homotopy class
of closed curves in $M$ or $N$ contains a unique closed geodesic. Hence
the marked length spectrum can be thought of as a map from conjugacy
classes in $\pi _1$ to the real numbers, assigning to each class the
length of the unique closed geodesic in that class. By having the
same marked length spectrum we now mean that there is an isomorphism
$\Psi :\pi_1(M) \to \pi_1(N)$ which pulls back the marked length spectrum
of $N$ to that of $M$.

Given an exact symplectomorphism $\phi$ of the open unit tangent
bundles, we claim that the induced map $\phi_* :\pi_1(M) \to \pi_1(N)$
gives such an isomorphism $\Psi$. To see this, observe again that
the closed geodesic $\gamma$ in a particular class $[\gamma]$ of
$\pi_1(M)$ can be lifted to a closed orbit of the characteristic flow
on the unit tangent bundle, since this characteristic flow is exactly
the geodesic flow. Now the symplectic homology theory, restricted
to considering orbits in the class $[\gamma]$, can be used to show that
the class $\phi_*[\gamma]$ similarly contains a unique closed orbit
of the characteristic flow of the same length. In other words, the
closed geodesic of $N$ in the class $\phi_*[\gamma]$ of $\pi_1(N)$
has the same length as the geodesic $\gamma$ in $M$. This establishes
our claim and proves the theorem.

\qed

\subsection{Relations with complex geometry}

\vskip 2mm

The authors first were drawn to this subject by the paper of
J.-C. Sikorav {\cite{sik}}, and comments on it made to us by David
Barrett concerning the translation of the symplectic rigidity results
there to holomorphic rigidity results. This leads directly to the
consideration of Grauert tubes \cite{gs}, \cite{ls}, \cite{bh}, a
natural complex structure on a ball bundle in the (co-)tangent bundle
of a Riemannian manifold. A translation of Theorem 1 above in terms of
complex structures is the following

\vskip 2mm

{\bf Corollary:}
For two compact Riemannian surfaces $X, X'$ of stricly negative curvature,
with Grauert tube complex structures defined on the respective tangent ball
bundles of radius $r$, $T_{r}(X), T_{r}(X')$, the following are equivalent:

\vskip 2mm

\hspace{10mm} (i) $T_r(X)$ and $T_r(X')$ are symplectomorphic 
	      
\vskip 1mm

\hspace{10mm} (ii) $T_r(X)$ and $T_r(X')$ are biholomorphic
              
\vskip 1mm

\hspace{10mm} (iii) $X$ and $X'$ are isometric.

\vskip 5mm

The result of Benci and Sikorav \cite{sik} gives a similar result, but
for translation invariant sets in $T(T^n)$, with fibers over $T^n$
which have vanishing first homolgy. The only 
Riemannian disk bundles in this case are for
flat metrics on $T^n$. Our rigidity result for surfaces holds when the
metric is allowed not to be of constant
sectional curvature, but the metric has to be
Riemannian. Much of the rest of the paper is directed at trying to
derive rigidity statements for non-Riemannian tubes in the tangent
bundle. We hope to return to the complex analytic aspects of the
non-Riemannian case at a future date. We wish to thank David Barrett
for his insightful remark.

\subsection{The case of genus $0$}

\vskip 2mm

We have not said anything in this paper about metrics on the $2$-sphere.
It turns out the Eliashberg and Hofer's construction in \cite{eh} can
be used to give an example of two arbitrarily small perturbations of the
round metric on the sphere such that the corresponding symplectic domains
in $T^* S^2$ have symplectomorphic interiors but non-conjugate Reeb flows
on the boundary, hence the closed domains are not symplectomorphic.

It is not clear however whether this construction can be performed such
that the perturbed metrics are Riemannian rather than just Finsler. It is
also unknown whether the round sphere itself is rigid, that is, if any
other Riemannian or Finsler metric has a symplectomorphic open unit
tangent bundle.  Symplectic homology will not provide the answer though.
Ideas from \cite{ziller} as applied in \cite{hs} can be used to give
examples of Riemannian metrics on $S^2$ whose unit tangent bundles have
the same symplectic homology (and volume) as the round metric.

\section{A technical lemma}

The purpose of this section is to prove the following lemma which
will be needed for our main result. Let $M$ be a $3$-manifold and
$X$ a never-vanishing vector field. We fix a metric on $M$ and
measure all norms with respect to it.

\begin{lemma}
There exists an $\epsilon$ (depending only on $M$, $X$ and the metric)
such that if $f$ is a diffeomorphism of $M$ with $\|f-\imath \|_{\infty}<\epsilon$ and $\|f_* X-X\|_{\infty}<\epsilon$ then $f$ is isotopic to the identity
through maps $f_t$ with $\|f_{t*} X-X\|_{\infty}<2\epsilon$.
\end{lemma}

Actually, using a theorem of J. Cerf, see \cite{cerf}, it can be
shown that any diffeomorphism of a $3$-manifold which is $C^0$-close
to the identity is in fact isotopic to the identity. However this is
a much deeper result than we need and anyway we will be using the
condition about the vector field.

\begin{proof}

We choose an open cover $B_1(x_i), 1\le i\le n$ for $M$ of 
coordinate balls of
radius $1$ and assume this is such that on each
ball we can choose coordinates $(x,y,z)$ with $X$ represented by
$\frac{\p}{\p z}$.  Further we assume that in fact the balls
$B_{\frac{1}{4}}(x_i)$ also cover $M$.

In each $B_1(x_i)$ we look at the disks $D^{\pm}_i =\{(x,y,\pm\frac{1}{2}|x^2+y^2 \le \frac{1}{2}\}$. We may assume that as $i$ varies these disks
are all disjoint in $M$. By taking $\epsilon$ sufficiently small, any given
$f$ can be isotoped to the identity in a neighbourhood of each $D^{\pm}_i$.
To do this, we note that $f_* X$ intersects each $D_i$ transversally so
we can project $f(D_i)$ along $f_* X$ and assume that $f$ maps each
$D_i$ to itself. The maps $f|_{D_i}$ must then be isotopic to the identity
and this isotopy can be extended to a small neighbourhood at the expense
perhaps of increasing $\|f_* X -X\|$ slightly.
The idea is now to isotope $f$ to the identity between the $D^{\pm}_i$.

\begin{lemma}
Suppose a disk in some $D_i$ is connected by flow-lines of $X$, say of length
no more than $1$, to a disk in $D_j$. Let $V$ be the union of the flow-lines.
Then there is an isotopy of $f$, which leaves $f$ equal to the identity near $D_i$ and $D_j$
and away from $V$, such that the resulting map is the identity on $V$ except
perhaps for an $\epsilon$-small neighbourhood of its boundary. Furthermore, if the
original $f$ leaves a subset of flow-lines in $V$ fixed, parameterized by a
subdomain of $D_i$ with finitely many components and a smooth boundary, then
the whole isotopy will leave these flow-lines unchanged.
\end{lemma}

\begin{proof}

This is an application of the following fact.

{\bf Fact} Given a diffeomorphism $g$ of the disk,
equal to the identity near the boundary,
there exists an isotopy of $g$ to the identity leaving the boundary fixed.
The construction can also be
carried out for a family $g_t, 0\le t \le 1$, that is, there exists a
corresponding smooth family of isotopies.
If the $g_t$ are also equal to the identity on a smooth
subdomain with finitely many components, then the family of isotopies
can be arranged also to leave these components fixed. 

Without the condition on leaving a certain subdomain of the disk
fixed, this is described exactly in \cite{thurston}, section $3.10$. 
To get the result as stated above one reduces to the
original case by first isotoping to the identity in a connected domain
enclosing the fixed-point set.

\begin{figure}
\centerline{\psfig{figure=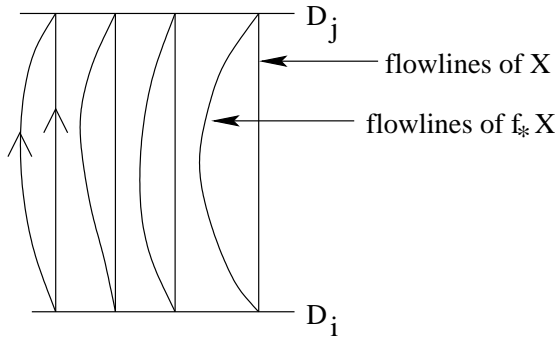}}
\caption{flowlines between two disks}
\end{figure}

The isotopy required by the lemma is now easily constructed. We start
with a homotopy between the vector fields $f_* X$ and $X$, say
$X_t$ where $X_0=f_* X$ and $X_1=X$.
This can be arranged so that $\|X_t-X\|_{\infty}<\epsilon$
for all $t$. We cut-off the $X_t$ to remain as $f_* X$ near the boundary of
$V$. Let $h_t$ be the diffeomorphism of a neighbourhood of $V$ defined
by leaving $D_i$ fixed and flowing along the vector field $-X$ to $D_i$ then back along $X_t$ for a similar time. Then $h_0=f$ and $h_1=\imath$ away from the
boundary of $V$, see figure $1$.
Unfortunately such $h_t$ are not
necessarily equal to the identity near $D_j$, except of course near its
boundary. However we can rescale $X_t$
to ensure that the map preserves $D_j$ and then compose with suitable isotopies
in the $(x,y)$-planes provided by the above fact to correct this.
On a subset of $V$ where $f_* X = X$ we may assume that $X_t =X$ for all
$t$ and the maps $h_t$ can be taken to be
the identity all along.

\qed
\end{proof}

We now apply this lemma repeatedly between different $D_i$ and $D_j$.
Notice that a map is equal to the identity if it is equal to the identity
on the flow-lines connecting sets of the form
$\tilde{D}^{\pm}_i =\{(x,y,\pm \frac{1}{2})|x^2+y^2\le\frac{1}{4}\}\subset D^{\pm}_i$ since these regions cover $M$. The order in which we apply the lemma must be chosen carefully
however, so as not to disturb regions in which $f$ has already been
isotoped to the identity.
We first isotope $f$ to the identity on the union $V_{ij}$ of flow-lines
between all $D_i$ and $D_j$ such that the flow-lines are of length less
than $1$ and the $V_{ij}$ do not intersect any other $D_k$. The disks
$D_k$ referred to here are the same as those $D^{\pm}_k$ defined above,
although we are not assuming that $D_i$ and $D_j$ are equal to $D^+_k$
and $D^-_k$ for the same $k$.

\begin{figure}
\centerline{\psfig{figure=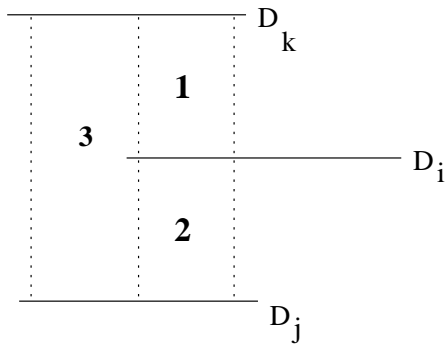}}
\caption{order of applying Lemma $6$}
\end{figure}

Figure $2$ divides the region between two disks $D_j$ and $D_k$ into
three numbered regions showing the order in which the isotopy provided
by Lemma $6$ should
be carried out if another disk $D_i$ intersects the flowlines between
them. 

Next we apply the lemma to all remaining pairs of $D^{\pm}_i$. The point
now is that if any $D_k$ happens to sit between $D^+_i$ and $D^-_i$ then,
away from an $\epsilon$-neighbourhood of its boundary, it must lie
on a complete set of flow-lines of $X$ from $D^+_i$ to $D^-_i$ on which $f$ has already
been isotoped to the identity. Thus the new isotopy will not affect $f$
here. After all of the above isotopies then, the resulting map
$f$ is equal to the identity.

\qed
\end{proof}
\section{Proof of main result}

\subsection{Construction of a diffeomorphism}

In this section we will construct a diffeomorphism between two closed
domains whose interiors are symplectomorphic.

We start with two domains $W_1$ and $W_2$ in $\mathcal M^{\circ}$. If
the interiors of $W_1$ and $W_2$ are symplectomorphic, then using
the symplectic homology theory and perhaps an application of Nielsen's
theorem as in the proof of Theorem $1$ we may assume that $\Sigma_1=\p W_1$ and
$\Sigma_2=\p W_2$ have the same marked length spectrum.
We observe that the differentials of diffeomorphisms of $M$ preserve the
set $\mathcal M^{\circ}$.

We choose a small $r$ such that $rW_1$ lies in the interior of both
$W_1$ and $W_2$. Then foliate the rest of $W_1$ by $s\Sigma_1$, 
$r\le s\le 1$, and the rest of $W_2$ by hypersurfaces $\Sigma_{2,s}$
where $\Sigma_{2,s}=s\Sigma_1$ for $s$ close to $r$ and
$\Sigma_{2,s}=s\Sigma_2$ for $s$ close to $1$. Further, all
$\Sigma_{2,s}\in \p \mathcal M$. This is possible for $r$ sufficiently small
by the connectedness
of $\mathcal M^{\circ}$ and the fact that all of our domains have
star-shaped fibers.

Now, as demonstrated originally by D. Anosov in \cite{anosov}, there
exists a continuous family of homeomorphisms $s\Sigma_1 \to \Sigma_{2,s}$
which map the Reeb flow (of the restriction of the Liouville form $\lambda$)
on $s\Sigma_1$ to that on $\Sigma_{2,s}$. Actually, as is made precise
in \cite{mor}, if we identify all of our foliating hypersurfaces with
a fixed $3$-manifold $\Sigma$, these homeomorphisms can be thought of as
the flow of a continuous time-dependent vector field, $Y_s$ say, on
$\Sigma$, where $Y_s$ can be differentiated only in the direction of
the Reeb flow at time $s$. We arrange things so that $Y_s \equiv 0$ for
$s$ close to $r$ and to $1$.

The homeomorphism $f$ for $s$ close to $1$ is a conjugacy between
Anosov flows which by assumption have the same marked length spectrum.
This homeomorphism is in fact H\"{o}lder continuous (see for
instance Chapter $9$ of \cite{kat}) and hence
we can apply a theorem of Livsic, see \cite{livsic}, which constructs a
function $g$ on $\Sigma$ such that the homeomorphism $\phi$ defined by
shifting a point $x$ a distance $g(x)$ along the flow-line through $x$
makes $\phi \circ f$ a time-preserving conjugacy, that is, it preserves
the Reeb vector field itself as opposed to just the flow-lines.

Such a conjugacy must in fact be of class $C^1$ by a result of
J. Feldman and D. Ornstein, see \cite{fo}, and so preserve the contact
form $\lambda$. We now use the theorem
of R. de la Llave and R. Moriyon in \cite{mor2} which says
that our time-preserving conjugacy must be smooth.
Now let $f$ be this diffeomorphism. The next step is to extend $f$,
thought of as a diffeomorphism of the boundaries, to a diffeomorphism
between the domains $W_1$ and $W_2$.

The homeomorphism $\phi$ is clearly the flow of a vector field, which
must be differentiable along the Reeb flow, so we may still assume that
$f$ is the time-$1$ map of a continuous vector field $Y_s$ still identically
zero for $s$ close to $1$.

On $s\Sigma_1$, $s$ close to $1$, we simply set $f(x)=sf(\frac{x}{s})$.

Away from $\Sigma_1$, we approximate $Y_s$ by a smooth vector field.
This can be done in such a way that the resulting one-parameter family
of diffeomorphisms, say $f_s$, $C^0$-approximate the original homeomorphisms
and map the Reeb flow on $s\Sigma_1$ to
flows $C^0$ close to the flows on $\Sigma_{2,s}$. Suppose that $f_s$ is
sufficiently close to $f$ on a level $1-2\delta$ that we can apply Lemma $6$.
Then the $f_s$ can be redefined for $1-2\delta \le s \le 1-\delta$ to be
the derived isotopy between $f_{1-2\delta}$ and $f$. These $f_s$ extend
$f$ smoothly over all levels $s\Sigma_1$ and $f$ further extends as
the identity inside $r\Sigma_1$.

\subsection{Isotopy of diffeomorphism into smooth symplectomorphism}

In this section we find an isotopy of $W_1$ which, composed
with the
$f$ of the previous section, gives a smooth symplectomorphism between $W_1$ and $W_2$.

\begin{lemma}

Suppose $\omega_0$ and $\omega_1$ are two symplectic (nondegenerate,
antisymmetric) bilinear forms on $\R^4$ and $\Sigma^3 \subset \R^4$
a linear subspace such that

$(i)$ if $v,w \in \Sigma$ and $\omega_0(v,w)$ or $\omega_1(v,w)$ is
positive, then $\omega_0(v,w)$ and $\omega_1(v,w)$ are nonnegative;

$(ii)$ if $v\in ker\omega_0|_{\Sigma}\cap ker\omega_1|_{\Sigma}$ and
$\omega_0(v,w)$ is positive, then $\omega_1(v,w)$ is nonnegative.

Then, for all $0\le t \le 1$, $t\omega_0 + (1-t)\omega_1$ is symplectic.
\end{lemma}

\begin{proof}

The kernel of any degenerate $t\omega_0 + (1-t)\omega_1$ must be at least
$2$-dimensional and so have non-trivial intersection with $\Sigma$,
which is clearly a contradiction.

\qed
\end{proof}

\begin{corollary}

For all $0\le t\le 1$, $t\omega + (1-t)f^*\omega$ is a symplectic
form on $W_1$, where $\omega=d\lambda$ is the canonical symplectic
form on the cotangent bundle.

\end{corollary}

\begin{proof}

For any $x\in W_1$, $f^*\omega$ and $\omega$ satisfy the conditions
of the above lemma on $T_x W_1$ with $\Sigma=T_x \Sigma_s$ as $f$ is
orientation preserving and approximately preserves the Reeb vector
field.
\qed

Now, since $f^*\lambda$ and $\lambda$ agree near $\p W_1$, we can apply
Moser's method to find an isotopy of $W_1$, fixed near $\p W_1$, which
generates a symplectomorphism between $f^*\omega$ and $\omega$.
Specifically, the isotopy can be taken to be the time-$1$ flow of the
time-dependent vector field $Z_t$ uniquely defined by
$Z_t\rfloor (t\omega + (1-t)f^*\omega)=\lambda - f^*\lambda$.
Note that $Z_t \equiv 0$ both near $\p W_1$ and near the zero-section.

The composition of this isotopy with our original diffeomorphism,
denoted again by $f$, is now the required symplectomorphism between
$W_1$ and $W_2$.

We now represent $f$ explicitly as a Hamiltonian diffeomorphism.
Observe that
associated to any $1$-form $\mu$ on $T^* M$ is, in the terminology of
\cite{elg}, a `contracting' vector field $X_{\mu}$ defined by 
$X_{\mu}\rfloor \omega =-\mu$. In the case of $\mu=\lambda$ or
$\mu=f^*\lambda$, this vector field vanishes only along the zero-section
$M$ and the associated flow contracts a disk towards each point on $M$.
For $X_{\lambda}$, these disks are just the cotangent fibers.
Now, $f$ maps $X_{f^*\lambda}$ into $X_{\lambda}$. The only map doing this
which is fixed near $M$ is defined as follows. Flow along $X_{f^*\lambda}$
until we are in the region where $f=\imath$, then flow out along
$-X_{\lambda}$ for the same time. 

Let $\phi_t$ and $\phi_t^{'}$ denote the
time-$t$ flows of $X_{\lambda}$ and $X_{f^*\lambda}$ respectively.
Assume that $\phi_T^{'}(W_1)$ lies in the region where $f=\imath$.
Note that as $f^*\lambda = \lambda$ near $\p W_1$ we can extend
$f^*\lambda$ and $X_{f^*\lambda}$ smoothly to $T^* M$.

Define an isotopy $h_t$, $0\le t\le T$ by $h_t=\phi_t^{-1}\circ \phi_t^{'}$.
Then $h_0=\imath$ and $h_T=f$.

Now, ${\cal L}_{X_{\lambda}}\omega =d(X_{\lambda}\rfloor \omega)=-\omega$
and similarly ${\cal L}_{X_{f^* \lambda}}\omega =-\omega$ so we have 
$\phi_t^*\omega=\phi_t^{'*}\omega =e^{-t}\omega$ and the $h_t$ are
all symplectomorphisms.

Let $V_t =\frac{dh_t}{dt}$, then $0=L_{V_t}\omega =d(V_t\rfloor \omega)$.
Hence the form $V_t \rfloor \omega$ is closed and the isotopy is
Hamiltonian if it is exact. But $V_t$ vanishes near the zero-section
and so we can use a parameterized version of the Relative Poincar\'{e}
Lemma to construct a smooth family of $H_t$ on $T^* M$ such that
$V_t \rfloor \omega =dH_t$ as required.

\end{proof}

\end{document}